\documentclass[12pt]{amsart}
\batchmode
\begin{document}
\numberwithin{equation}{section}

\def\1#1{\overline{#1}}
\def\2#1{\widetilde{#1}}
\def\3#1{\widehat{#1}}
\def\4#1{\mathbb{#1}}
\def\5#1{\frak{#1}}
\def\6#1{{\mathcal{#1}}}

\def\C{{\4C}}
\def\R{{\4R}}
\def\N{{\4N}}
\def\Z{{\4Z}}

\title[]{A Burns-Krantz type theorem for domains with corners}
\author[L.~Baracco, D.~Zaitsev, G.~Zampieri]{Luca~Baracco, Dmitri~Zaitsev, and Giuseppe~Zampieri}
\footnotetext{Dipartimento di Matematica, Universit\`a di Padova, via Belzoni 7, 35131 Padova, Italy,
{\tt baracco@math.unipd.it, zaitsev@math.unipd.it, zampieri@math.unipd.it}}
%\subjclass{}

\maketitle
%\tableofcontents

%\def\Label#1{\label{#1}{\bf (#1)}~}
\def\Label#1{\label{#1}}

% Standard sets

\def\cn{{\C^n}}
\def\cnn{{\C^{n'}}}
\def\ocn{\2{\C^n}}
\def\ocnn{\2{\C^{n'}}}

% Abbreviations

\def\const{{\rm const}}
\def\rk{{\rm rank\,}}
\def\id{{\sf id}}
\def\aut{{\sf aut}}
\def\Aut{{\sf Aut}}
\def\CR{{\rm CR}}
\def\GL{{\sf GL}}
\def\Re{{\sf Re}\,}
\def\Im{{\sf Im}\,}

\def\codim{{\rm codim}}
\def\crd{\dim_{{\rm CR}}}
\def\crc{{\rm codim_{CR}}}

\def\phi{\varphi}
\def\eps{\varepsilon}
\def\d{\partial}
\def\a{\alpha}
\def\b{\beta}
\def\g{\gamma}
\def\G{\Gamma}
\def\D{\Delta}
\def\Om{\Omega}
\def\k{\kappa}
\def\l{\lambda}
\def\L{\Lambda}
\def\z{{\bar z}}
\def\w{{\bar w}}
\def\Z{{\1Z}}
\def\t{\tau}
\def\th{\theta}

\emergencystretch15pt
\frenchspacing

\newtheorem{Thm}{Theorem}[section]
\newtheorem{Cor}[Thm]{Corollary}
\newtheorem{Pro}[Thm]{Proposition}
\newtheorem{Lem}[Thm]{Lemma}

\theoremstyle{definition}\newtheorem{Def}[Thm]{Definition}

\theoremstyle{remark}
\newtheorem{Rem}[Thm]{Remark}
\newtheorem{Exa}[Thm]{Example}
\newtheorem{Exs}[Thm]{Examples}

\def\bl{\begin{Lem}}
\def\el{\end{Lem}}
\def\bp{\begin{Pro}}
\def\ep{\end{Pro}}
\def\bt{\begin{Thm}}
\def\et{\end{Thm}}
\def\bc{\begin{Cor}}
\def\ec{\end{Cor}}
\def\bd{\begin{Def}}
\def\ed{\end{Def}}
\def\br{\begin{Rem}}
\def\er{\end{Rem}}
\def\be{\begin{Exa}}
\def\ee{\end{Exa}}
\def\bpf{\begin{proof}}
\def\epf{\end{proof}}
\def\ben{\begin{enumerate}}
\def\een{\end{enumerate}}

\section{Introduction}

In their pioneering paper \cite{BK94}, {\sc Burns} and {\sc Krantz} have established a boundary version
of the classical {\sc Cartan}'s uniqueness theorem: 

\medskip
{\em If $f$ is holomorphic self-map of a smoothly bounded strongly pseudoconvex domain $D\subset\C^N$ 
such that $f(z)=z+o(|z-p|^3)$ as $z\to p$ for some $p\in\d D$, then $f(z)\equiv z$.}
\medskip

Furthermore, they give an example showing that the exponent $3$
in the above statement cannot be decreased.
More recent results of this kind are due to {\sc Vlacci-Tauraso} \cite{VT98} and
{\sc Bracci-Vlacci-Tauraso} \cite{BVT00} in one and to {\sc Huang} \cite{H95} and {\sc Gentili-Migliorini} \cite{GM97}
in several complex variables. A basic global assumption in these results was that
the given holomorphic map is a self-map of the given bounded domain.

The goal of this paper is twofold. 
First, to give purely local boundary uniqueness results for maps defined only 
on one side as germs at a boundary point and hence not necessarily sending any domain to itself
and also under the weaker assumption that $f(z)=z+o(|z-p|^3)$ holds only
for $z$ in a proper cone in $D$ with vertex $p$. 
Such results have no analogues in one complex variable in contrast to the 
situation when a domain is preserved.
And second, to extend the above results from boundaries of domains to submanifolds of higher codimension.
Here the usual replacement for a {\em one-sided neighborhood} of $M$ is a {\em wedge with the edge} $M$.
The notion of strong pseudoconvexity can be also extended to this situation (see \S\ref{preli} for more details).
We prove:

\bt\Label{mainhigh}
Let $M\subset\C^N$, $N\ge 2$, be a generic submanifold at a point $p\in M$
and $f$ be a germ at $p$ of a holomorphic self-map of a strongly pseudoconvex wedge with edge $M$ at $p$
such that 
%\begin{equation}\Label{asymptotics}
$$f(z)=z+o(|z-p|^3)$$
%\end{equation} 
as $z$ approaches $p$ nontangentially. 
Then $f(z)\equiv z$.
Moreover, we can restrict the way $z$ approaches $p$
by requiring $z$ to belong to a cone
with vertex $p$,
only depending only on the domain of definition of $f$,
which is properly contained in the given wedge.
\et

Note that the map $f$ is not assumed to have any limit at points of the edge $M$ other than $p$
(though it has nontangential limits almost everywhere by a result of {\sc Forstneri\v c} \cite{Fo92})
and, in case it has a continuous boundary value, it is not assumed to send $M$ into itself.
The authors of this paper are not aware of any other rigidity results for maps of this kind.
Here the known method based on the {\sc Lempert} theory of extremal discs 
for the Kobayashi metric \cite{L81a,L81b,L82} cannot be applied.
Instead, our approach makes use of the theory of stationary discs
recently developed by {\sc Tumanov} \cite{T01}.

We conclude by giving a consequence of Theorem~\ref{mainhigh} in the hypersurface case
which is a ``local'' and ``nontangential'' version of corresponding results
of \cite{BK94,H95}:

\bc\Label{mainhyp}
Let $M\subset\C^N$, $N\ge 2$, be a strongly pseudoconvex hypersurface at a point $p\in M$
and $f$ be a germ at $p$ of a holomorphic self-map of the pseudoconvex side
such that $f(z)=z+o(|z-p|^3)$ as $z$ approaches $p$ nontangentially. 
Then $f(z)\equiv z$.
\ec

In fact we shall obtain Theorem~\ref{mainhigh} and Corollary~\ref{mainhyp}
as special cases of the more general Theorem~\ref{main} below,
where the source and the target of $f$ may be different even as germs at $p$.

\section{Preliminaries and generalizations}\Label{preli}
Let $M\subset\C^N$ be a smooth submanifold.
Recall that $M$ is called {\em generic} in $\C^N$
if its tangent space at each point spans $\C^N$ over complex numbers.
Denote by $T^*\C^N$ the {\em cotangent bundle} of $\C^N$
(regarded as the space of all $(1,0)$ forms on $\C^N$)
and, for every $p\in M$, by $N_p M:=T_p\C^N/T_pM$ the {\em normal} and by
$N^*_p M\subset T^*_p \C^N$ 
the {\em conormal} spaces to $M$ in $\C^N$,
where $N^*_p M$ consists of all $1$-forms $\phi\in T^*_p\C^N$
such that $\Re \phi|_{T_pM} =0$. Recall that the
(vector-valued) {\em Levi form}
of $M$ at $p\in M$, $L_p\colon T^{1,0}_p M\times T^{1,0}_p M\to N_p M\otimes \C$
is the (unique) hermitian form such that the equality
$$L_p(X_p,Y_p)=\frac{1}{2i}[X,\1Y]_p \mod T^{1,0}_pM\oplus T^{0,1}_pM$$
holds for all vector fields $X$ and $Y$ in $T^{1,0}M$,
where $X_p$ stands for the evaluation of the vector field $X$ at $p$.
The {\em Levi cone} $C_p=C_p M\subset N_p M$ is the convex hull
of all vectors $L_p(v,v)$ for $v\in T^{1,0}_pM\setminus \{0\}$
and the {\em dual Levi cone} is $C_p^*M:=\{\xi\in N_p^*M : \xi|_{C_p M}>0\}$.
The notion of strong pseudoconvexity can be extended to
submanifolds of higher codimension as follows.
The submanifold $M\subset\C^N$ is {\em strongly pseudoconvex at $p\in M$} 
(see e.g. \cite{TH83,Fo91,T01}) if $C_p^*M\ne\emptyset$,
i.e.\ if $\xi\big(L(u,v)\big)$ is positive definite for some conormal $\xi$.
It follows that $M$ is strongly pseudoconvex at $p$ if and only if 
there exists a strongly pseudoconvex (real) hypersurface in $\C^N$
that contains a neighborhood of $p$ in $M$ (see e.g.\ \cite{Fo91}).

We now introduce the notion of a strongly pseudoconvex wedge with edge $M$ at $p$
which generalizes strongly pseudoconvex sides of hypersurfaces
and which is relevant for the rigidity result given by Theorem~\ref{main}.
For open cones $\G,\G'\subset\R^n$, we write $\G\ll\G'$ if 
$\1{\G}\setminus\{0\}\subset \G'$ and say that {\em $\G$ is properly contained in $\G'$}.
More generally, a cone $\G\subset\R^n$ is said to be {\em proper} in an open subset $D\subset\R^n$
if it is properly contained in another cone $\G'\subset\R^n$ such that a neighborhood of $0$ 
in $\G'$ is contained in $D$.
By a {\em wedge with edge $M$ at $p$} in the direction of an open cone $\G\subset N_p M$
we shall understand here any domain $W\subset\C^N$ such that
for all open cones $\G',\G''\subset N_p M$ with $\G' \ll \G \ll \G''$
and for every sufficiently small neighborhood $U$ of $p$ in $\C^N$, one has:
\begin{gather}
(M\cap U) + (\G'\cap U) \subset W;\label{sub} \\
(M\cap U) + (\G''\cap U) \text{ contains a neighborhood of $p$ in $W$},\label{sup}
\end{gather}
where the normal space $N_pM$ is identified with any fixed complementary subspace to $T_pM$ in $\C^N$.
It is easy to see that the given definition is independent of the choice of such a complement
and also of the choice of local coordinates.
Conditions \eqref{sub} and \eqref{sup} give estimates on the shape of $W$ from inside and from outside
respectively and clearly depend only on the intersection of $W$ with an arbitrarily small neighborhood of $p$
in $\C^N$.

We call $\G$ the {\em directional cone} of $W$ at $p$.
We say that the wedge $W$ is {\em strongly pseudoconvex} at $p$ if its directional cone $\G$ satisfies the following:
\begin{gather}
L_p(v,v)\in \G \text{ for some vector } v\in T^{1,0}_pM;\label{substr}\\
\xi \text{ is positive on } \1\G\setminus\{0\} \text{ for some covector } \xi \in C_p^* M.\label{supstr}
\end{gather}
Here $L_p$ denotes the Levi form of $M$ at $p$ as above.
As in the above definition of a wedge, also here \eqref{substr} and \eqref{supstr}
are estimates on the shape of $\Gamma$ from inside and from outside respectively.
In particular, it follows from \eqref{supstr} that $C_p^* M$ is nonempty 
and hence $M$ is automatically strongly pseudoconvex in the above sense.
An important class of examples can be obtained as follows.

\begin{Exa}
Let $\rho_1,\ldots,\rho_d$, $1\le d\le N-1$, 
be strictly plurisubharmonic functions in a neighborhood of $p$ in $\C^N$,
vanishing at $p$ and satisfying $\d\rho_1\wedge\cdots\wedge \d\rho_d\ne 0$.
Then 
$$W:=\{\rho_1<0,\ldots,\rho_d<0\}$$ 
is a strongly pseudoconvex wedge at $p$ with edge 
$M:=\{\rho_1=\ldots=\rho_d=0\}$.
In particular, an one-sided neighborhood of a real hypersurface
is a strongly pseudoconvex wedge if and only if it is the strongly pseudoconvex
side in the usual sense.
\end{Exa}

Another important class is given by Siegel domains of 2nd kind (see \cite{Pi})
that are strongly pseudoconvex wedges whose edges are their Shilov boundaries.

The following is a generalization of Theorem~\ref{mainhigh}:

\bt\Label{main}
Let $M\subset\C^N$, $N\ge 2$, be a generic submanifold through a point $p$,
$U$ and $V$ be strongly pseudoconvex wedges with edge $M$ at $p$
and $f$ be a germ at $p$ of a holomorphic map between $U$ and $V$ with $f(z)=z+o(|z-p|^3)$
as $z$ approaches $p$ nontangentially.
Then $f(z)\equiv z$.
\et

It will follow from the proof that one has the same statement under the weaker assumptions
that $U$ satisfies only \eqref{sub} and \eqref{substr} (for some cone $\G\subset N_pM$) 
and $V$ satisfies only \eqref{sup} and \eqref{supstr} (for a possibly different cone $\2\G\subset N_pM$).
Also, we can weaken the assumptions about the asymptotics of $f$
by asking that $f(z)=z+o(|z-p|^2)$ and $\langle \xi,f(z)\rangle = \langle \xi, z\rangle + o(|z-p|^3)$
for at least one covector $\xi$ such that \eqref{supstr} holds.
We shall write $U$ (resp. $V$) for wedges with edge $M$ in the statements below
that will be used for the source (resp. target) wedge in Theorem~\ref{main}.

\section{One-dimensional case}
One of the key points in the proofs of 
mentioned rigidity results in \cite{BK94} and \cite{H95}
is the reduction to the one-dimensional situation,
where $D=\D$ is the unit disc in $\C$.
The latter case is therefore of particular importance.
We begin with a new elementary proof of the following result 
due to \cite{BK94,H95}:

\bp\Label{elem}
Let $f\colon\D\to \D$ be a holomorphic self-map with $f(z)=z + o(|z-1|^3)$
as $z\to 1$ in $\D$. Then $f(z)\equiv z$.
\ep

\bpf
By the classical Fatou's theorem, 
$f$ has an $L^\infty$ boundary value function on $\d\D$ that we denote also by $f$.
Since $f(z)\in\1\D$ for all $z\in \1\D$, we have 
\begin{equation}\Label{7a0}
\Re \Big( \bar z\frac{z-f(z)}{|z-1|^{4}} \Big) \geq 0, \quad z\in\d\D.
\end{equation}
We write $K_\eps(1)$ for the disc with center $1$ and radius $\eps>0$.
Then
\begin{multline*}%\Label{8a}
\Re  \int_{\{\theta: e^{i\theta}\notin K_\eps(1) \} }
 e^{-i\theta} \frac{e^{i\theta} -f(e^{i\theta})}{|e^{i\theta}-1|^{4}} d\theta
=\Im \int_{ \d\Delta\setminus K_\eps(1) }
z^{2}{\bar z} \frac{z-f(z)}{(z-1)^{4} } \frac {dz}{z} \\
=\Im \int_{ \d K_\eps(1) \cap \Delta }
z \frac{z-f(z)}{(z-1)^{4} } \frac {dz}{z} \to 0, \quad \eps\to 0,
\end{multline*}
where the second equality holds because the function under the integral 
extends holomorphically to $\Delta$
and the convergence to $0$ is a consequence of the estimate $f(z)=z+o(|z-1|^3)$.
Hence we must have in \eqref{7a0} the equality almost everywhere on $\d\Delta$
proving $f(z)\equiv z$ as required.
\epf

The next result, for which we give a self-contained proof,
generalizes Proposition~\ref{elem}.

\bp\Label{best}
Let $f\colon\D\to \D$ be a holomorphic self-map, 
$z_k\in \D$ a sequence converging nontangentially to $1$
such that $f(z_k)=z_k + o(|z_k-1|^{3})$ as $k\to\infty$.
Then $f(z)\equiv z$.
\ep

\bpf
We begin by deriving a special case of the classical Julia's Theorem following \cite{R80}.
For $a\in \D$, consider the automorphism $\phi_a(z):=(a-z)/(1-\bar az)$ of $\D$
interchanging $a$ and $0$. Then we have 
\begin{equation}\Label{schwarz}
|\phi_{f(a)}(f(z))|\le |\phi_a(z)|
\end{equation}
for all $a,z\in \D$ as a consequence of the Schwarz lemma.
We next use the identity
\begin{equation}\Label{id}
1-|\phi_a(z)|^2 = \frac{(1-|a|^2)(1-|z|^2)}{|1-\bar az|^2}
\end{equation}
that can be verified by easy computation (see \cite[Theorem~2.2.2]{R80}).
Then (\ref{schwarz}--\ref{id}) imply
\begin{equation}\Label{}
\frac{|1-\1{f(a)} f(z)|^2}{1-|f(z)|^2} \le \frac{1-|f(a)|^2}{1-|a|^2}\, \frac{|1-\bar az|^2}{1-|z|^2}.
\end{equation}
Setting $a=z_k$ and taking the limit we obtain the estimate
\begin{equation}\Label{julia}
\frac{|1-f(z)|^2}{1-|f(z)|^2} \le \frac{|1-z|^2}{1-|z|^2}
\end{equation}
which is the desired special case of Julia's Theorem.

The second part of the proof closely follows the proof of Lemma~2.1 in \cite{H95}.
An easy computation shows that the harmonic function
\begin{equation}\Label{}
\xi(z):=\Re \Big(\frac{1+z}{1-z} - \frac{1+f(z)}{1-f(z)} \Big)
\end{equation}
is nonnegative in view of \eqref{julia}. Furthermore, we have $|(f(z_k)-z_k)/(1-z_k)|<1$ for $k$ sufficiently large
and hence, using the geometric series expansion,
\begin{equation}\Label{}
\frac{1+f(z_k)}{1-f(z_k)}=\frac{(1+z_k+f(z_k)-z_k)/(1-z_k)}{1-(f(z_k)-z_k)/(1-z_k)}=
\frac{1+z_k}{1-z_k} + o(|z_k-1|), 
\quad k\to\infty.
\end{equation}
Hence $\xi(z_k)=o(|z_k-1|)$ as $k\to\infty$. 
In view of the Hopf lemma, the latter fact is only possible if $\xi(z)\equiv 0$
and therefore $f(z)\equiv z$ as required.
\epf

\section{Stationary discs}
We review basic fact from the {\sc Tumanov}'s theory of stationary discs in higher 
codimension extending the classical {\sc Lempert}'s theory for strictly (linearly) convex domains.
Recall that a continuous mapping $\Phi\colon\1\Delta\to \C^N$,
where $\Delta$ is the unit disc in $\C$, is called
an {\em analytic disc attached to $M$}
if $\Phi|_{\Delta}$ is holomorphic and $\Phi(\d\Delta)\subset M$.
In this paper we shall consider only smooth analytic discs.
In his celebrated paper \cite{L81a}, {\sc Lempert} 
gave a characterization of those analytic discs attached to a boundary
of a strictly convex domain, for which the maximum of the Kobayashi distance is attained,
in terms of meromorphic lifts of its boundary values. 
The discs admitting such lifts were called {\em stationary}.
More recently, {\sc Tumanov} \cite{T01} extended this notion to higher codimension.
In his terminology, an analytic disc $\Phi$ attached to $M$ is {\em stationary} if there exists
a nonzero holomorphic map (a lift) $\Phi^*\colon\1\Delta\setminus\{0\} \to T^*\C^N$
such that $\zeta \Phi^*(\zeta) \in \6O(\Delta)\cap \6C(\1\Delta)$
and $\Phi^*(\zeta)\in N^*_{\Phi(\zeta)} M$ for $\zeta\in\d\Delta$.
In contrast to the full space of attached analytic discs,
the space of stationary discs (satisfying certain inequalities) 
has finite dimension provided $M$ is sufficiently nondegenerate
(see Theorem~\ref{tumanov-stationary} below).

In order to have an explicit description of stationary discs,
it will be convenient to write $M$ near a point $p\in M$ in the form 
\begin{equation}\Label{form}
M = \{\rho(x+iy,w)= 0\} = \{ y=h(x,w)\}, \quad \rho(x+iy,w)= h(x,w)-y,
\end{equation}
where coordinates $(z,w)=(x+iy,w)\in\C^d\times\C^n=\C^N$ vanishing at $p$ are chosen
such that $h(0)=0$ and $h'(0)=0$. 
In order to have a description of all holomorphic (and meromorphic) lifts of an attached 
analytic disc $\Phi$, it will be convenient, following \cite{T88,T01}, to introduce
``partial holomorphic lifts'' of $\Phi$ to $T^*\C^N$ (or their collections forming a matrix)
for which only the ``$z$-components'' are holomorphic. 
More precisely, to every analytic disc $\Phi$ sufficiently small (e.g.\ in $\6C^\infty(\1\Delta)$)
which is attached to $M$ (given in the form \eqref{form}),
one can associate, by solving a Bishop's equation, 
a unique smooth real invertible $d\times d$ matrix function $G(\zeta)$ on $\d\Delta$
such that $G(1)=\id$ and $G(\zeta)\rho_z(\Phi(\zeta))$ extends holomorphically from $\d\D$ to $\D$,
where $\rho_z:=(\rho_{z_1},\ldots,\rho_{z_d})$ is the gradient with respect to the chosen coordinates.
(More precisely, the function $G(\zeta)$ here is normalized at $\zeta=1$ rather than $\zeta=0$
as in  \cite{T01} and, hence, differs from it by an invertible constant matrix factor.)
We now state {\sc Tumanov}'s existence and uniqueness result for stationary discs
in the following slightly adapted for our purposes form:

\bt[Tumanov]\Label{tumanov-stationary}
Let $M\subset\C^N$ be a smooth generic submanifold through a point $p$.
Then for every fixed $c^0\in \R^d$
such that the Levi form of $M$ at $p$ is nondegenerate in the direction of
the conormal $c^0\d\rho$,
every fixed $\eps>0$ and every data $(\l,c,v)\in \C^d\times\R^d\times \C^n$
sufficiently close to $(0,c^0,0)$,
there exists a unique smooth stationary disc $\Phi(\zeta)=(z(\zeta),w(\zeta))$, $\zeta\in\1\Delta$,
with $\Phi(1)=p$ and $w'(1)=v$ whose lift $\Phi^*$
is of the form 
\begin{equation}\Label{lift-form}
\Phi^*(\zeta)=\big(\Re(\lambda\zeta+c)\big)G(\zeta)\d\rho(\Phi(\zeta))
\end{equation}
and such that 
$w(\zeta)=\a(\zeta)+\b(\zeta)$ with $\a$ linear and $\|\b\|_{\6C^\infty}<\eps\|v\|$.
\et

Since Theorem~\ref{tumanov-stationary} is not explicitly stated in this form in \cite{T01},
we provide here a proof, essentially following the arguments of \cite[\S4--5]{T01}.

\bpf
Let $M$ be given near $p=0$ by \eqref{form} and
let $\Phi\colon \1\D\to\C^N$ be a sufficiently small stationary disc for $M$ with a lift 
$\Phi^*\colon\1\Delta\setminus\{0\} \to T^*\C^N$.
Denote by $G(\zeta)$ the invertible $d\times d$ matrix function on $\d\D$ as described above
such that $G(1)=\id$ and $G(\zeta)\rho_z(\Phi(\zeta))$ extends holomorphically to $\Delta$.
Observe that $G(\zeta)$ is arbitrarily close to the constant identity provided $\Phi$ is arbitrarily small.
We can then write 
\begin{equation}\Label{try}
\Phi^*(\zeta)=\nu(\zeta)\,G(\zeta)\,\big(\rho_z(\Phi(\zeta)),\rho_w(\Phi(\zeta))\big), \quad \zeta\in\d\D,
\end{equation}
where $\nu\colon \d\D\to \R^d$ is a real vector function.
By definition of the lift, $\zeta\Phi^*(\zeta)$ extends holomorphically to $\D$.
Recall that the construction of $G$ required $\Phi$ to be small.
We can also further require that 
the inverse $S(\zeta):=(G(\zeta)\rho_z(\Phi(\zeta)))^{-1}$ exists
and extends holomorphically to $\D$. Multiplying $\zeta\Phi^*(\zeta)$ by $S(\zeta)$ on the right
and using \eqref{try} we conclude that the product $\zeta \nu(\zeta)$ must also extend holomorphically to $\D$.
Since $\nu$ is valued in $\R^d$, expanding into the Fourier series on $\d\D$, 
we conclude that $\nu(\zeta)=\Re(\l\zeta+c)$
for suitable vectors $\l\in \C^d$, $c\in \R^d$.
In particular, $\Phi^*(\zeta)$ is always of the form \eqref{lift-form}.

We now explore, under what conditions, an attached analytic disc $\Phi$
with a lift $\Phi^*$ given by \eqref{lift-form} is stationary, i.e.\ when does
$\zeta\Phi^*(\zeta)$ extend holomorphically to $\D$.
Since
\begin{equation}\Label{extend}
2\zeta \Re(\l\zeta+c) = \zeta (\l\zeta + \bar\l \bar\zeta + 2c) = \l\zeta^2 + \bar\l + 2c\zeta,
\quad \zeta\in\d\D,
\end{equation}
obviously has a holomorphic extension to $\D$ and since
$G(\zeta)\rho_z(\Phi(\zeta))$ extends holomorphically to $\D$  by the construction, then
\begin{equation}\Label{}
\zeta\big(\Re(\lambda\zeta+c)\big)G(\zeta)\rho_z(\Phi(\zeta)),
\quad \zeta\in\d\D,
\end{equation}
also has the same extension property. Hence it remains to find $\Phi(\zeta)$ such that
the ``$w$-component'' 
\begin{equation}\Label{w-comp}
\zeta\big(\Re(\lambda\zeta+c)\big)G(\zeta)\rho_w(\Phi(\zeta)),
\quad \zeta\in\d\D,
\end{equation}
has the extension property. The fact that $\Phi(\zeta)=(z(\zeta),w(\zeta))$ is attached to $M$ with $\Phi(1)=0$
can be expressed via the Bishop's equation
\begin{equation}\Label{bishop}
\Re z(\zeta) = -T_1 \big( h(\Re z(\zeta), w(\zeta))\big),
\quad \zeta\in\d\D,
\end{equation}
so that $w(\zeta)$ with $w(1)=0$ can be seen as an independent parameter
and $z(\zeta)$ is uniquely determined from it by \eqref{bishop} and the condition $z(1)=0$.

We next assume, without loss of generality, that $h(x,w)$ is further normalized as 
\begin{equation}\Label{}
h(x,w)=(w^* A_j w)_j + o(|(x,w)|^2), \quad (x,w)\to 0,
\end{equation}
where Hermitian $n\times n$ matrices $A_1,\ldots,A_d$ represent the ($\R^d$-valued) Levi form of $M$ at $0$,
and that the corresponding normalization holds also for the derivatives:
\begin{equation}\Label{rho}
\rho_w(z,w) = h_w (z,w) = (\bar A_j \bar w)_j + o(|(x,w)|), \quad (x,w)\to 0,
\end{equation}
where we think of $L_w$ as a $d\times n$ matrix
(we write $*$ for  the conjugate transposed matrix).
The next step consists of substituting \eqref{rho} into \eqref{w-comp}
and writing $G(\zeta)=\id+G_0(\zeta)$ and
\begin{equation}\Label{matrices}
\sum_{j=1}^d\big(\Re(\l_j\zeta+c_j)\big) A_j =
 (\id -\zeta X)^* B (\id-\zeta X),
\quad\zeta\in\d\D,
\end{equation}
for suitable $n\times n$ complex matrices $B=B(\l,c)$, $X=X(\l,c)$.
To show the existence of $B$ and $X$, we set 
$P=P(\l,c):=\sum_j \l_j A_j$, $Q=Q(\l,c):=\sum_j c_j A_j$
and rewrite \eqref{matrices} as
\begin{equation}\Label{}
\zeta P + \bar\zeta P^* + 2Q = B - \zeta BX - \bar\zeta X^*B + X^*BX
\end{equation}
leading to the system
\begin{equation}\Label{matrix-eqs}
 P^*X^2 + 2QX + P = 0, \quad 2Q=B+X^*BX.
\end{equation}
Recall that we consider $(\l,c)$ in a neighborhood of $(0,c^0)$ such that 
$P(0,c^0)=0$ and $Q(0,c^0)$ is invertible.
Then the system \eqref{matrix-eqs} can be solved for $(X,B)$ near $(0,B^0)$ by the implicit function theorem
with $B^0=2Q(0,c^0)$ being invertible.

We now look for $w(\zeta)$ in the form 
\begin{equation}\Label{look}
w(\zeta)=w_0+\zeta (\id-\zeta X)^{-1} \sigma(\zeta),
\end{equation}
where $\sigma(\zeta)$ is an unknown function
and $w_0$ is uniquely determined by the condition $w(1)=0$.
Note that, if $(X,B)$ is sufficiently close to $(0,B^0)$, the matrix
$(\id-\zeta X)$ is invertible for $\zeta\in\1\D$.
Substituting \eqref{look} 
(and the corresponding solution $z(\zeta)$ of \eqref{bishop})
into \eqref{rho} and using \eqref{matrices} and the property of $G(\zeta)$, 
we conclude that the expression in \eqref{w-comp} can be written as
\begin{equation}\Label{}
 \zeta\big(\Re(\lambda\zeta+c)\big) (\bar A_j \bar w_0)_j 
+ (\id-\bar\zeta \bar X)^* \bar B \1{\sigma(\zeta)}
+ o(\|\sigma\|),
\quad \zeta\in\d\D,
\end{equation}
where the norm of $\sigma$ in the last term is understood
in any H\"older space $\6C^\a$, $0<\a<1$.
The first term always extends holomorphically to $\D$ in view of \eqref{extend}.
Hence the holomorphic extendibility of \eqref{w-comp} is equivalent to that of
\begin{equation}\Label{red}
(\id-\bar\zeta \bar X)^* \bar B \1{\sigma(\zeta)}
+ o(\|\sigma\|),
\quad \zeta\in\d\D.
\end{equation}
Furthermore, since the matrix factor in \eqref{red}
is invertible and holomorphic in $\zeta$, our condition is equivalent
to the antiholomorphic extendibility of $\sigma(\zeta) + o(\|\sigma\|)$.
By taking a Cauchy transform $K$ given by the Cauchy integral formula
up to a constant
and using the fact that $\sigma$ is holomorphically extendible,
we obtain an equation $\sigma(\zeta)=a+K(o(\|\sigma\|))$,
where $a\in \C^n$ is a parameter that can be chosen to be $\sigma(1)$.
The latter equation can be solved by the implicit function theorem
in any H\"older class $\6C^{k,\a}$ for any sufficiently small $a=\sigma(1)$
and any data $(\l,c)$ sufficiently close to $(0,c^0)$.
It can be now seen from the construction and the implicit function theorem
that the solution $w(\zeta)$
given by \eqref{look} defines a smooth transformation 
$(\l,c,a)\mapsto (\l,c,w'(1))$ sending $(0,c^0,0)$ into itself
and being locally invertible there.
This shows the existence and uniqueness of a small solution $w(\zeta)$
with given data $\l$, $c$ and $v=w'(1)$ having the required properties.
The proof is complete.
\epf

Note that, by the construction in \cite{T01},
for $c$ fixed and $\l,v$ arbitrarily small,
the difference $\Phi^*(\zeta)-\xi$ is also arbitrarily small on $\d\Delta$,
where $\xi:=\Phi^*(1)$. We shall say that a stationary disc $\Phi$
has a {\em small lift $\Phi^*$ in the direction of a conormal} $\xi\in N^*_pM$
if $\Phi(1)=p$, $\Phi^*(1)=\xi$ and both $\Phi(\zeta)-p$, $\Phi^*(\zeta)-\xi$
are small in the norm $\6C^\infty(\d\Delta)$.
Then the uniqueness part of Theorem~\ref{tumanov-stationary} implies that, 
if $\xi\in C^*_p$, all stationary discs with small lifts in the direction $\xi$ are
given by Theorem~\ref{tumanov-stationary} with sufficiently small data $\l,v$.

We now show that the assumption \eqref{substr} in \S\ref{preli}  for a wedge with edge $M$
implies that the corresponding stationary discs are contained in the wedge.

\bl\Label{contain}
Let $U$ be a wedge with edge $M$ at $p$ whose directional cone is $\G$
and $v_0\in T_p^{1,0}M$ be such that $L_p(v_0,v_0)\in \G$.
Then, for any $\xi\in C_p^*$, any neighborhood of the origin
contains an open subset of parameters $\l,v$ 
for which the corresponding stationary discs $\Phi$ with small lifts
in the direction $\xi$ satisfy $\Phi(\Delta)\subset U$.
Moreover, for any sector $S\subset\Delta$ with vertex $1$
and sufficiently small parameters $\l,v$, the images $\Phi(S)$ 
are contained in a fixed proper cone in $U$
with vertex $p$.
\el

\bpf
By Theorem~\ref{tumanov-stationary}, there exist stationary discs 
$\Phi(\zeta)=(z(\zeta),w(\zeta))$ with arbitrarily small lifts in the direction $\xi$ with 
$w'(1)\in\C v_0$. Since also $w(\zeta)=\a(\zeta)+\b(\zeta)$ with $\a$ and $\b$
satisfying the conclusion of Theorem~\ref{tumanov-stationary}
and since $\Phi$ is attached to $M$, the derivatives $\Phi'(\zeta)$
will be arbitrarily close to a multiple of $v_0$ for all $\zeta\in\1\Delta$.
Choose any conormal $c\d\rho(p)\in N^*_p M$ with $c\, d\rho(p;L_p(v_0,v_0))>0$, where $c\in\R^d$
and $\rho$ is the defining function of $M$ near $p$.
Then, for $\Phi$ as above, the function $c\rho\circ\Phi$ is plurisubharmonic in $\Delta$
and is zero on the boundary $\d\Delta$. By the maximum principle, $c\rho(\Phi(\zeta))\le 0$
for all $\zeta\in\Delta$. By choosing finitely many conormals $c_1\d\rho,\ldots,c_s\d\rho$
such that 
$$L_p(v_0,v_0)\in\bigcap_j\{c_j\, d\rho(p)>0\}\subset\subset \G,$$ 
we conclude that $\Phi(\Delta)\subset \cap_j \{c_j\rho<0\}\subset U$
for $\Phi$ sufficiently small. 
The last conclusion follows from the Hopf lemma
applied to suitable small perturbations of $\rho$.
\epf

\bl\Label{identity}
Let $\Phi_0$ be a stationary disc with sufficiently small lift in the direction $\xi\in N^*_p M$.
Then, for stationary discs $\Phi$ arbitrarily close to $\Phi_0$ in the direction $\xi$,
the union of their boundaries $\Phi(\d\Delta)$ contains open subsets in $M$ arbitrarily close to $p$.
\el

\bpf
By Theorem~\ref{tumanov-stationary}, the discs $\Phi$
are in one-to-one correspondence with the parameters $\l\in\C^d$, $v\in\C^n$.
We then consider the correspondence
$$
(\l,v)\underset\psi\mapsto(\Phi'(1),(\Phi^{*})'(1))\underset\pi\mapsto\Phi'(1),
\quad \C^d\times\C^n\to T_{(p,\xi)} N^*M \to T_pM.
$$
Since $M$ is strongly pseudoconvex, $\Phi$ has defect $0$,
i.e.\ its lift $\Phi^*$ cannot be chosen holomorphic.
By \cite[Proposition~3.9]{T01}, $\psi$ is injective
and hence immersive on a dense subset.
Since $\dim_\R N^* M=2(n+d)$, $\psi$ is also submersive on a dense subset
of its domain of definition. Then $\pi\circ\psi$ is clearly submersive on the same subset
and the required statement follows.
\epf

Lemmas~\ref{contain} and \ref{identity} imply:

\bc\Label{cor-id}
In the setting of Lemma~\ref{contain} the union of images $\Phi(\Delta)$ of the 
stationary discs $\Phi$ in the direction $\xi$ arbitrarily close to $\Phi_0$,
contains smooth generic submanifolds of $\C^N$ arbitrarily close to $p$.
\ec

\section{Strongly convex hypersurfaces passing through stationary discs}

Following {\sc Lempert} \cite{L81a} we call an analytic disc $\Phi\colon\1\Delta\to\C^N$
stationary for a bounded domain with smooth boundary $\Omega\subset\C^N$
if it is stationary for $M=\d\Omega$ and the lift $\Phi^*$ can be chosen such that
it does not vanish on $\d\Delta$.
Our goal here is to construct, for each sufficiently small stationary disc $\Phi$,
a smooth strongly convex domain for which $\Phi$ is also stationary with the same lift.

\bp\Label{pass}
Suppose that a wedge $V$ with edge $M$ at $p$ and a conormal $\xi\in C_p^*$ satisfies \eqref{supstr}.
Then there exist a biholomorphic change of coordinates in $\C^N$,
a neighborhood $B(p)$ of $p$ in $\C^N$ and, for every stationary disc $\Phi$ with 
sufficiently small data $\l,v$ as in Theorem~{\rm\ref{tumanov-stationary}}
and $\Phi^*(1)$ close to $\xi$,
a strongly convex (with respect to the new coordinates) bounded domain $\Omega\subset\C^N$ such that
\begin{enumerate}
\item[(i)] $V\cap B(p) \subset \Omega$;
\item[(ii)] $M\cap B(p) \subset \d\Omega$;
\item[(iii)] $\Phi^*(\zeta)|_{T_{\Phi(\zeta)}\d\Omega}=0$ and $\Phi^*(\zeta)\ne 0$ for all $\zeta\in\d\Delta$.
\end{enumerate}
In particular, $\Phi$ is stationary also for $\Omega$ with the same lift $\Phi^*$.
\ep

\bpf
In the vector notation we have, for every stationary disc $\Phi$ as above,
$\Phi^*(\zeta)= c(\zeta)\partial\rho(\Phi(\zeta))$
for $\zeta\in\d\Delta$, where $c(\zeta)$ is a real vector and
$\rho$ is the defining function of $M$ as in \eqref{form}.
If the data $\l,v$ are sufficiently small, it follows from Theorem~\ref{tumanov-stationary}
that $\Phi\colon\1\Delta\to \C^N$ is a smooth embedding.
Hence we can find a smooth extension 
$\2c$ of $c$ from $\Phi(\d\Delta)\cong \d\Delta$ to a neighborhood of $p$ in $M$.
The neighborhood can be chosen uniformly for all $\Phi$ as above. 
Moreover, since $\Phi^*(\zeta)$ is arbitrarily close to $\xi$
and $\xi$ satisfies \eqref{supstr},
also the extension $\2c \d\rho$ can be chosen to have these properties.
Then one can see that the real hypersurface defined by $\2\rho:=\2c \rho + C\|\rho\|^2$
for a sufficiently large constant $C>0$ near $p$ can be extended to the boundary
of a domain $\Omega$ as required.
\epf

In the situation of Proposition~\ref{pass}, we can take advantage of 
{\sc Lempert}'s theory \cite{L81a,L81b},
in particular, of the fact that stationary discs coicide with geodesics
for the Kobayashi metric and also that the latter are regular up to the boundary.

\section{A direct proof of Theorem~\ref{main} in a weaker form}

Here we give a direct proof of the statement of Theorem~\ref{main},
where the asymptotics $f(z)=z + o(|z-p|^{3})$ is 
assumed for any $z$ in the wedge
rather than in a proper cone.
We begin by a uniqueness result for stationary discs in smooth strongly convex domains.

\bp\Label{special}
Let $\Omega\subset\C^N$ be a smooth strongly convex bounded domain 
and  $\Phi\colon\1\Delta\to \1\Omega$ be a stationary disc for $\Omega$.
Let $\2\Phi\colon \Delta\to \Omega$ be a holomorphic map
satisfying $\2\Phi(\zeta)=\Phi(\zeta) + o(|\zeta-1|^3)$ as $\zeta\to 1$ in $\Delta$. Then  $\2\Phi\equiv \Phi$.
\ep

\bpf
Since $\Phi$ is stationary for $\Omega$,
there exists a (smooth) lift $\Phi^*\colon \1\Delta\setminus\{0\} \to T^*\C^N$
such that $\Phi^*(\zeta)\in N^*_{\Phi(\zeta)}\d\Omega\setminus\{0\}$ for $\zeta\in \d\Delta$
and $\zeta\Phi^*(\zeta)$ is holomorphic in $\Delta$.
By the classical Fatou's theorem, $\2\Phi$ has an $L^\infty$ boundary value on $\d\Delta$ that we also denote by $\2\Phi$.
%We can assume $\2\Phi(\zeta)\in\1\Omega$ for all $\zeta\in\d\Delta$.
Since $\Omega$ is strongly convex, we can choose $\Phi^*$ such that
\begin{equation}\Label{7a}
\Re\Big\langle\Phi^{*}(\zeta), \frac{\Phi(\zeta) -\2\Phi(\zeta)}{|\zeta -1|^{4}} \Big\rangle \geq 0
\end{equation}
for almost all $\zeta\in\d\Delta$,
where $\langle\cdot,\cdot\rangle$ denotes the standard pairing.
Furthermore, since $\Omega$ is strongly convex, for almost every $\zeta\in\partial\Delta$, we have the equality in \eqref{7a} 
if and only if $\Phi(\zeta)=\2\Phi(\zeta)$. Since $\Phi$ is smooth up to the boundary of $\Delta$
and because of the given estimate for $\2\Phi$,
we can repeat the arguments of the proof of Proposition~\ref{elem} to obtain
\begin{equation*}%\Label{8a}
 \Re \int_{\{\theta : e^{i\theta}\notin K_\eps(1)\}}
\Big\langle\Phi^{*}(e^{i\theta}), \frac{\Phi(e^{i\theta}) -\2\Phi(e^{i\theta})}{|e^{i\theta}-1|^{4}} \Big\rangle d\theta
\to 0, \quad \eps\to 0.
\end{equation*}
Hence we must have the equality in \eqref{7a} almost everywhere on $\d\Delta$
which proves $\Phi \equiv \2\Phi$ on $\Delta$ as required.
\epf

We now turn to the proof of Theorem~\ref{main} in the weaker form mentioned above:

\bp\Label{mainweak}
The conclusion of Theorem~{\rm\ref{main}} holds 
if the hypothesis on the asymptotics of $f$ at $p$
holds without restrictions on the way of approaching $p$.
\ep

\bpf
By Lemma~\ref{contain}, for any $\xi\in C^*_p$,
there exists an open set of stationary discs $\Phi$
with arbitrarily small lifts in the direction $\xi$
such that $\Phi(\Delta)\subset U$.
In view of Corollary~\ref{cor-id},
it is enough to show that $f\circ \Phi \equiv \Phi$ on $\Delta$
for each disc $\Phi$ as above.
For $\Phi$ fixed, let $B(p)$ and $\Omega$ be given by Proposition~\ref{pass}.
Since $f(z)\to p$ as $z\to p$ and since $\Phi$ can be taken arbitrarily small,
we may assume that $\Phi(\Delta)\subset V\cap B(p)\subset\Omega$.
Since $\Phi$ is smooth up to the boundary, we have 
the estimate $f(\Phi(\zeta))=\Phi(\zeta)+o(|\zeta-1|^3)$
and hence we are in the situation of Proposition~\ref{special}
with $\2\Phi:=f\circ\Phi$ that yields $\2\Phi\equiv\Phi$ as required.
\epf

\section{Uniqueness for images of stationary discs}

In order to prove Theorem~\ref{main} as it stands 
we shall need results of {\sc Lempert} \cite{L81a,L81b,L82} for stationary discs and geodesics 
and a result of {\sc Huang} \cite{H95} on boundary uniqueness for holomorphic self-maps of the disc $\Delta$.
The main consequence can be stated as follows.

\bp\Label{unique}
Let $\Omega$ and  $\Phi$ be as in Proposition~{\rm\ref{special}}
and $\2\Phi\colon \Delta\to \Omega$ be a holomorphic map
satisfying $\2\Phi(z_k)=\Phi(z_k)+o(|z_k-1|^3)$ 
for some sequence $z_k$ in $\Delta$ converging to $1$ nontangentially.
Then $\2\Phi \equiv \Phi$.
\ep

\bpf
Since $\Phi$ is stationary for $\d\Omega$,
there exists a holomorphic retraction $\pi\colon \Omega\to\Delta$
such that $\pi\circ\Phi=\id_\Delta$ (see \cite{L81a,L81b,L82})
and $\pi$ is smooth up to the boundary of $\Omega$.
Hence $\pi\circ\2\Phi$ is a holomorphic self-map of $\Delta$
satisfying $\pi(\2\Phi(z_k))=z_k+o(|z_k-1|^3)$.
Then $\pi \circ \2\Phi = \id$ by Proposition~\ref{best}. 
Since both $\pi$ and $\2\Phi$
are holomorphic, we have
\begin{equation}\Label{kobayashi}
K_\Delta\big(\pi(\2\Phi(\zeta_1)),\pi(\2\Phi(\zeta_2))\big)
\le K_\Omega(\2\Phi(\zeta_1),\2\Phi(\zeta_2)) \le K_\Delta(\zeta_1,\zeta_2),
\end{equation}
where $K_D$ denotes the Kobayashi distance in a domain $D$.
Since $\pi \circ \2\Phi = \id$, the first and the last distances in \eqref{kobayashi}
are the same. We conclude that $K_\Omega(\2\Phi(\zeta_1),\2\Phi(\zeta_2)) = K_\Delta(\zeta_1,\zeta_2)$
holds for all $\zeta_1,\zeta_2\in\Delta$ and hence $\2\Phi$ is a complex geodesic of $\Omega$.
Since $\Omega$ is strongly convex and smooth, it follows from \cite{L81a,L81b}
that $\2\Phi$ is smooth up to the boundary of $\Delta$.
The regularity of $\Phi$ and $\2\Phi$ on $\1\Delta$ 
together with Taylor's formula yields the improved estimate 
$\2\Phi(\zeta)=\Phi(\zeta)+o(|\zeta-1|^3)$ for $\zeta\to p$
without any restriction on the way of approaching $p$.
The required conclusion follows from Proposition~\ref{special}.
Alternatively, one can use at this last stage uniqueness results for complex geodesics.
\epf

\section{Proofs of Theorems~\ref{mainhyp}--\ref{main}}

The proof of Theorem~\ref{main} is obtained 
by repeating the proof of Proposition~\ref{mainweak}
and using Proposition~\ref{unique} instead of Proposition~\ref{special}.
Note that we enter the hypotheses of Proposition~\ref{unique}
in view of the last conclusion of Lemma~\ref{contain}.
Theorem~\ref{mainhigh} is a special case of Theorem~\ref{main} with $U=V$.

\end{document}